\documentclass{ifacconf}

\usepackage{graphicx}      
\usepackage{natbib}        
\usepackage{tikz} 
\usepackage{url} 
\usepackage{amsmath} 
\usepackage{amssymb} 
\usepackage[implicit=false]{hyperref}

\usepackage{algorithm}
\usepackage{enumerate}
\usepackage[noend]{algpseudocode}
\usepackage[utf8]{inputenc}
\usepackage[caption=false, font=footnotesize]{subfig}
\usetikzlibrary{shapes,arrows,calc}
\usetikzlibrary{calc}
\usetikzlibrary{angles, quotes}
\usetikzlibrary{decorations.pathmorphing,decorations.markings,patterns}
\usepackage{todonotes}
\usepackage{pgfplots}
\usepgfplotslibrary{fillbetween}

\newcommand{\mydefinition}{defn}
\newcommand{\myremark}{rem}
\newcommand{\mylemma}{lem}
\newcommand{\myproof}{pf}

\DeclareMathOperator*{\argmin}{arg\,min}
\DeclareMathOperator*{\sol}{sol}
\begin{document}
\definecolor{set19c1}{HTML}{E41A1C}
\definecolor{set19c2}{HTML}{377EB8}
\definecolor{set19c3}{HTML}{4DAF4A}
\definecolor{set19c4}{HTML}{984EA3}
\definecolor{set19c5}{HTML}{FF7F00}
\definecolor{set19c6}{HTML}{FFFF33}
\definecolor{set19c7}{HTML}{A65628}
\definecolor{set19c8}{HTML}{F781BF}
\definecolor{set19c9}{HTML}{999999}

\pgfplotstableread{data/random_times.dat}{\randomTimes}
\pgfplotstableread{data/random_iters.dat}{\randomIters}
\pgfplotstableread{data/random_nodes.dat}{\randomNodes}
\pgfplotstableread{data/result_invpend.dat}{\invpend}
\begin{frontmatter}

\title{BnB-DAQP: A Mixed-Integer QP Solver for Embedded Applications} 

\thanks[footnoteinfo]{This work was supported by the Swedish Research Council (VR) under contract number 2017-04710.}

\author[First]{Daniel Arnström}
\author[First]{Daniel Axehill}

\address[First]{Division of Automatic Control, Link\"oping University, Linköping, Sweden (e-mail: daniel.\{arnstrom,axehill\}@liu.se).}

\begin{abstract} 
We propose a mixed-integer quadratic programming (QP) solver that is suitable for use in embedded applications, for example, hybrid model predictive control (MPC). The solver is based on the branch-and-bound method, and uses a recently proposed dual active-set solver for solving the resulting QP relaxations. Moreover, we tailor the search of the branch-and-bound tree to be suitable for embedded applications on limited hardware; we show, for example, how a node in the branch-and-bound tree can be represented by only two integers. The embeddability of the solver is shown by successfully running MPC of an inverted pendulum on a cart with contact forces on an MCU with limited memory and computing power. 
\end{abstract}


\end{frontmatter}
\section{Introduction}

In model predictive control (MPC), an optimization problem is solved at each time step to determine an optimal control action. When MPC is used to control safety-critical systems in real time, the employed optimization solvers need to be reliable and efficient, demands that become particularly challenging when considering embedded systems due to limited computational resources and memory. 
For MPC of \textit{linear} systems with \textit{continuous} states and controls, the optimization problems in question are commonly convex quadratic programs (QPs), for which there exist several reliable and efficient solvers that have been developed specifically for real-time MPC \citep[e.g.][]{patrinos2013accelerated,ferreau2014qpoases,frison2020hpipm,arnstrom2022daqp}. In \textit{hybrid} MPC, where some states and/or controls are restricted to take binary values, the resulting optimization problems are instead \textit{mixed-integer} QPs (MIQPs). Since MIQPs are nonconvex they cannot be solved as efficiently and reliably as their continuous/convex counterpart. 

A popular framework for solving MIQPs is that of \textit{branch and bound} (B\&B) \citep{land1960automatic}, where several QP relaxations of the nominal MIQP are solved in sequence. What differentiate B\&B solvers is how they determine which QP relaxations to solve, and how these are solved. 
A survey of strategies for selecting relaxations to solve is given in \cite{achterberg2005branching}. In the context of MPC, methods have been proposed that solve the relaxation using active-set methods \citep{axehill2006mixed,bemporad2015solving,bemporad2018numerically, hespanhol2019structure}, gradient projections methods \citep{axehill2008dual,naik2017embedded}, operator splitting methods \citep{stellato2018embedded}, and interior-point methods \citep{frick2015embedded,liang2020early}.

Although B\&B often finds an optimal solution sufficiently fast in practice for medium-sized problems, solving MIQPs is NP-hard, i.e., the worst-case complexity is limiting. This theoretical worst-case complexity, hence, often restricts B\&B methods from being used in real-time applications. As an alternative to B\&B methods,  heuristic approaches based on, for example, ADMM \citep{takapoui2020simple} and machine learning \citep{bertsimas2022online} have been proposed to solve MIQPs.  
While these methods show impressive computational time and often gives a sufficient, albeit suboptimal, solution, there is no formal guarantees on the solution quality, which make the resulting control law unreliable and, hence, unsuitable for MPC of safety-critical systems. 

Another emerging research direction that addresses the (theoretical) conservative worst-case complexity of B\&B methods is to develop complexity certification methods tailored for MPC. Such methods determine upper bounds on the worst-case number of computations required to solve any MIQP encountered in a \textit{given} hybrid MPC application \citep{axehill2010improved,shoja2022overall}. 

In this paper we present an MIQP solver that is based on B\&B and that uses the dual active-set solver DAQP \citep{arnstrom2022daqp} for solving relaxations, combined with a search strategy, described in Section \ref{ssec:branch}, that is tailored for embedded applications (by prioritizing simplicity). 
The proposed method falls directly into the complexity framework proposed in \cite{shoja2022overall}, enabling an overall worst-case certification of the solvers complexity given an MPC application.

The solver exploits the well-known warm-starting capabilities of active-set methods when solving sequences of similar QP relaxations in B\&B. Moreover, considering a \textit{dual} active-set method, in contrast to, for example, the \textit{primal} method used in \cite{hespanhol2019structure}, has two advantages. First, the dual solution to a relaxation is always a feasible starting point for a subsequent relaxation, while a primal solution is not (assuming that a binary constraint is fixed after solving a relaxation, see Section \ref{ssec:branch} for details). Secondly, a dual method allows for early termination when solving relaxations \citep{fletcher1998numerical}, which often saves a lot of computational effort. 
Compared with the active-set methods in \cite{bemporad2015solving,bemporad2018numerically}, which also can be interpreted as \textit{dual} active-set methods, the proposed method avoids some overhead stemming from a nonnegative least-squares reformulation used therein; see Sec. III.A in \cite{arnstrom2022daqp} and Remark \ref{rem:comp-nnls} herein for details.

The main contributions of the paper are: (i) An open-source C implementation of an MIQP solver for embedded applications, available under a permissive license and with support for complexity certification;
(ii) Use of a least-distance formulation of the relaxations to reduce computations (for example when computing upper bounds); (iii) A compact representation of the search tree; only $3 n_b$ integers are maximally needed to represent the tree, where $n_b$ is the number of binary constraints; (iv) A simple, yet effective, way of regularizing the Hessian when binary variables do not nominally enter the objective function (common in hybrid MPC, where such ``auxiliary'' binary variables originate from logical rules and switches).  

\section{Preliminaries}
\subsection{Problem formulation}
We consider problems of the form
\begin{subequations}
  \label{eq:miqp}
  \begin{align}
   &\underset{x}{\text{minimize}}&&\frac{1}{2} x^T H x + f^T x \label{eq:miqp-obj}\\
   &\text{subject to} &&\underline{b} \leq A x \leq \bar{b} \label{eq:miqp-con}\\
   & && A_i x \in \left\{ \underline{b}_i, \bar{b}_i  \right\},\quad \forall i \in \mathcal{B}\label{eq:miqp-bin}, 
  \end{align}
\end{subequations}
with decision variable $x\in \mathbb{R}^n$. The objective function \eqref{eq:miqp-obj} is characterized by $H\in \mathbb{S}^{n}_{++}$ and $f\in \mathbb{R}^n$ (how to handle a singular $H$ is addressed in Section \ref{sssec:ex-reg}); the feasible set \eqref{eq:miqp-con} is characterized by $A \in \mathbb{R}^{m\times n}$ and $\overline{b},\underline{b} \in \mathbb{R}^m$. The binary constraints \eqref{eq:miqp-bin}, which make \eqref{eq:miqp} a nonconvex problem, are given by $\mathcal{B} \subseteq \{1,\dots, m\}$ with $|\mathcal{B}| = n_b$. For a vector $v$ (matrix $M$), we denote $v_i$ ($M_i$) its $i$:th element (row). 

Specifically, \eqref{eq:miqp} generalizes mixed-integer quadratic programming problems, where a subset of the decision variables $x$ are binary, i.e., when $x_i \in \{0,1\}$ for $i\in \mathcal{B} \subseteq \{1,\dots,n\}$. Such problems are encountered in, for example, hybrid MPC of mixed logical dynamical (MLD) systems \citep[see, e.g.,][]{bemporad1999control}.

Instead of solving \eqref{eq:miqp} directly, we transform it into a least-distance problem (LDP) of the form

\begin{equation}
 \begin{aligned}
   \label{eq:mildp}
   &\underset{u}{\text{minimize}}&&\frac{1}{2} \|u\|_2^2 \\
   &\text{subject to} &&\underline{d} \leq M u \leq \bar{d} \\
   & && M_i u \in \left\{\underline{d}_i, \bar{d}_i \right\},\quad \forall i \in \mathcal{B}, 
  \end{aligned} 
\end{equation}
by doing the coordinate transformation $u = R x + v$, where $R$ is an upper triangular Cholesky factor of $H$ {(i.e., $H = R^T R$)} and $v\triangleq R^{-T} f$. Consequently, $M$, $\underline{d}$ and $\bar{d}$ are given by.
\begin{equation}
  \label{eq:aux-def}
  M \triangleq A R^{-1}, \quad \underline{d}\triangleq \underline{b}+M v, \quad \bar{d}\triangleq \bar{b}+M v.
\end{equation}

This transformation reduces intermediary computations in the proposed solver, presented in Section \ref{sec:main}, in particular when computing the upper bounds described in Section \ref{sssec:early-term}. Moreover, the up-front cost for transforming \eqref{eq:miqp} into \eqref{eq:mildp} is negligible for MIQPs encountered in hybrid MPC applications, which are often small to medium-sized, compared with the required computations for solving the relaxations described in Section \ref{ssec:solve-relax} (and solving relaxations of the transformed problem is often cheaper than solving relaxations of \eqref{eq:miqp}). Finally, for \textit{linear} hybrid MPC problems, the transformation can be done \textit{a priori} since $H$ and $A$ remain constant, resulting in no additional overhead online due to the transformation.

\begin{\myremark}
    \label{rem:comp-nnls}
      The relaxations solved in \cite{bemporad2015solving,bemporad2018numerically} can also be interpreted as LDPs, but these methods require the LDP solution to be transformed back into ``normal'' coordinates every time a relaxation is solved, resulting in a significant overhead (see Sec. III.A in \cite{arnstrom2022daqp} for details). Our approach only transforms the final, global, solution back into normal coordinates, i.e., it only performs a single coordinate transformation.
\end{\myremark}
\subsection{Branch and Bound}
A naive approach for solving \eqref{eq:mildp} would be to enumerate all possible combinations arising from the binary constraints in $\mathcal{B}$ and select the feasible solution among these subproblems with the smallest norm. Such a brute-force approach do, however, require $2^{|\mathcal{B}|}$ LDPs to be solved, which quickly becomes intractable when the number of binary constraints increases. 

In branch-and-bound (B\&B) methods \citep{land1960automatic}, all of these combinations are \textit{implicitly} considered by gradually fixing the constraints in $\mathcal{B}$, encoded by the sets $\underline{\mathcal{B}} \subseteq \mathcal{B}$ and $\overline{\mathcal{B}}\subseteq \mathcal{B}$ that contain indices corresponding to binary constraints that have been fixed at its lower and upper bound, respectively. After such fixations, LDP relaxations of the form  
\begin{equation}
 \begin{aligned}
   \label{eq:ldp-relax}
   &\underset{u}{\text{minimize}}&&\frac{1}{2} \|u\|_2^2 \\
   &\text{subject to} &&\underline{d} \leq M u \leq \bar{d} \\
   & && M_i u = \underline{d}_i,\quad \forall i \in \underline{\mathcal{B}} \\
   & && M_i u = \bar{d}_i,\quad \forall i \in \overline{\mathcal{B}} 
  \end{aligned} 
\end{equation}
are solved, and their solutions are used to dismiss other binary combinations that cannot (based on Lemma \ref{lem:dom} below) be optimal, which avoids explicitly solving the corresponding relaxations. To systematically consider such binary combinations, relaxations of the form \eqref{eq:ldp-relax} can be ordered in a tree, where a node in this tree is defined in the following way: 
\begin{\mydefinition}[Node]
    A \textit{node} is a pair $(\underline{\mathcal{B}},\overline{\mathcal{B}})$ with the sets $\underline{\mathcal{B}}, \overline{\mathcal{B}}\subseteq \mathcal{B}$ satisfying $\underline{\mathcal{B}} \cap \overline{\mathcal{B}} = \emptyset$. The \textit{level} of a node is given by a mapping $\ell: \mathbb{P}(\mathcal{B}) \times \mathbb{P}(\mathcal{B}) \to \mathbb{Z}_{\geq 0}$ defined by the rule $\ell(\underline{\mathcal{B}},\overline{\mathcal{B}})\triangleq  |\underline{\mathcal{B}}|+|\overline{\mathcal{B}}|$, where $\mathbb{P}(\mathcal{B})$ is the power set of $\mathcal{B}$.
\end{\mydefinition}

By \textit{processing}, or \textit{exploring}, a node $(\underline{\mathcal{B}}, \overline{\mathcal{B}})$ we mean solving the corresponding LDP relaxation in \eqref{eq:ldp-relax}. The ``gradual fixing'' of binary constraints mentioned above corresponds to moving down the tree, which corresponds to processing \textit{descendants} to nodes that have already been processed.
\begin{\mydefinition}[Descendant]
    A node $(\underline{\mathcal{B}}_d,\overline{\mathcal{B}}_d)$ is said to be a \textit{descendant} to the node $(\underline{\mathcal{B}},\overline{\mathcal{B}})$ if $\underline{\mathcal{B}} \subseteq \underline{\mathcal{B}}_d$, $\overline{\mathcal{B}} \subseteq \overline{\mathcal{B}}_d$, 
    and $\ell(\underline{\mathcal{B}_d},\overline{\mathcal{B}}_d) > \ell(\underline{\mathcal{B}},\overline{\mathcal{B}})$. 
    Moreover $(\underline{\mathcal{B}}_d,\overline{\mathcal{B}}_d)$ is a \textit{child} to $(\underline{\mathcal{B}},\overline{\mathcal{B}})$ (and, conversely, $(\underline{\mathcal{B}}, \overline{\mathcal{B}})$ is a \textit{parent} to $(\underline{\mathcal{B}}_d,\overline{\mathcal{B}}_d)$) if 
$\ell(\underline{\mathcal{B}_d},\overline{\mathcal{B}}_d)- \ell(\underline{\mathcal{B}},\overline{\mathcal{B}})=1$.
\end{\mydefinition}

Making the tree exploration more concrete, after processing a node $(\underline{\mathcal{B}},\overline{\mathcal{B}})$ in a B\&B method, an index $i\in \mathcal{B}: i\notin \underline{\mathcal{B}}$ and $i\notin \overline{\mathcal{B}}$ is selected and the descendants $(\underline{\mathcal{B}}\cup\{i\},\overline{\mathcal{B}})$ and $(\underline{\mathcal{B}},\overline{\mathcal{B}}\cup\{i\})$ are added to a list, denoted $\mathcal{T}$, which contain \textit{pending} nodes, i.e., nodes that are to be processed. This is the \textit{branching} step in B\&B. Selecting which index to branch over, and in which order pending nodes are processed, are described in more detail in Section~\ref{ssec:branch}.

To avoid explicitly enumerating all possible nodes, the solutions to previous relaxations can sometimes be used to dismiss the branching, commonly known as \textit{cuts} in the tree, while ensuring that an optimal solution is still obtained. These cuts are based on the following, well-known, lemma:

\begin{\mylemma}[Dominance]
  \label{lem:dom}
  Let $J$ be the optimal objective function value of a relaxation corresponding to a node $(\underline{\mathcal{B}},\overline{\mathcal{B}})$, and let $J_d$ be the optimal objective function value of one of its descendants; then $J \leq J_d$. 
\end{\mylemma}
\begin{\myproof}
  Directly follows from the feasible set of a child being a subset of the feasible set of its parent's feasible set (since more equality constraints are enforced in descendants). \qed
\end{\myproof}
Specifically, Lemma \ref{lem:dom} gives rise to two types of cuts: binary feasibility cuts and dominance cuts.
\subsubsection{Binary feasibility}
If a solution $\tilde{u}$ to a relaxation satisfies the binary constraints, i.e., if $M_i \tilde{u} \in \{\underline{d}_i,\bar{d}_i\}$, $\forall i\in \mathcal{B}$, no further descendant need to be explored since Lemma~\ref{lem:dom} implies that doing so could only lead to worse binary feasible solutions. 
\subsubsection{Dominance}
If $\bar{J}$ is the objective function value of the best binary feasible solution found so far, and the solution to a relaxation for a particular node yields $J \geq \bar{J}$, no further descendant need to be explored since, again, Lemma \ref{lem:dom} implies that solving descendant nodes can only lead to worse binary feasible solutions (infeasibility is a special case of this if we use the convention of $J=\infty$ for infeasible problems). Dominance cuts can also sometimes be invoked before a relaxation is solved completely if a dual ascent method is used as inner solver (see Section \ref{sssec:early-term} for details.) 

\subsubsection{General branch-and-bound method}
Based on the above concepts (solving relaxations, branching, and cuts), a generic formulation of a branch-and-bound method for solving \eqref{eq:mildp} is provided in Algorithm \ref{alg:bnb}, where relaxations are solved at Step \ref{step:solve-relax}, and the dominance/feasibility cuts are invoked at Step \ref{step:dom-cut} and \ref{step:feas-cut}. 
Still, several steps in Algorithm \ref{alg:bnb} remain to be made concrete; namely, how the relaxations at Step \ref{step:solve-relax} are solved, and in what order the tree is explored (determined by the selections at Steps \ref{step:node-selection}, \ref{step:branch-selection}, \ref{step:child-selection}.) 
Particular implementations, suitable for embedded applications, of Steps \ref{step:node-selection}, \ref{step:solve-relax}, \ref{step:branch-selection}, and \ref{step:child-selection} are the subject of the next section. 

\begin{algorithm}
  \caption{Generic B\&B method for solving \eqref{eq:mildp}}
  \label{alg:bnb}
  \begin{algorithmic}[1]
	\Require $M,\underline{d},\bar{d}, \mathcal{B}$ 
    \Ensure $u^*, J^*$ 
	\State $\bar{u}\leftarrow \star$, $\bar{J} \leftarrow \infty$, $\mathcal{T} \leftarrow \{(\emptyset, \emptyset)$\}
	\While{$\mathcal{T}\neq \emptyset$}
    \State $(\underline{\mathcal{B}},\overline{\mathcal{B}}) \leftarrow$ select node from $\mathcal{T}$ \label{step:node-selection}
    \State $u, J \leftarrow$ solve \eqref{eq:ldp-relax}\label{step:solve-relax}
    \If{$J \geq  \bar{J}$} \textbf{continue} \label{step:dom-cut}
    \Comment dominance
	\EndIf
    \If{$M_i u \in \{\underline{d}_i,\bar{d}_i\}, \forall i \in \mathcal{B}$} \label{step:feas-cut}
	\State $\bar{u} \leftarrow u, \bar{J} \leftarrow J$
    \Comment binary feasibility
	\Else
    \State $i\leftarrow$ select $i\in \mathcal{B}: M_i u \notin \{\underline{d}_i,\bar{d}_i\}$ \label{step:branch-selection}
    \State add $(\underline{\mathcal{B}}\cup\{i\},\overline{\mathcal{B}})$ and $(\underline{\mathcal{B}},\overline{\mathcal{B}}\cup\{i\})$ to $\mathcal{T}$ \label{step:child-selection}
	\EndIf
	\EndWhile
	\State \Return $u^* \leftarrow \bar{u}$, $J^* \leftarrow J$ 
  \end{algorithmic}
\end{algorithm}

\section{Proposed MIQP solver}
\label{sec:main}
In this section we outline the proposed MIQP solver, which is the main contribution of this paper, by making the steps in the generic B\&B method in Algorithm \ref{alg:bnb} specific. In particular, we describe how the relaxations at Step \ref{step:solve-relax}  are solved (Section~\ref{ssec:branch}), and how the selections at Steps \ref{step:node-selection}, \ref{step:branch-selection} and \ref{step:child-selection} are made in Algorithm \ref{alg:bnb} (Section~\ref{ssec:branch}).
All of these specifications coalesce into Algorithm \ref{alg:bnb-full}, presented in Section \ref{ssec:complete-method}.
\subsection{Solving relaxations}
\label{ssec:solve-relax}
To solve the LDPs in \eqref{eq:ldp-relax} we use the dual active-set solver DAQP \citep{arnstrom2022daqp} given in Algorithm~\ref{alg:daqp}, where $\lambda$ denotes the dual iterate. Here, we highlight some aspects that are essential for incorporating it in a branch-and-bound method; for a detailed description we refer the reader to \cite{arnstrom2022daqp}.
In particular, early-termination and double-sided/equality constraints were mentioned in \cite{arnstrom2022daqp},  but were never presented jointly in a single algorithm (which we do here in Algorithm \ref{alg:daqp}.)

\begin{algorithm}
  \caption{The dual active-set method from \cite{arnstrom2022daqp}, extended to handle equality constraints, double-sided constraints, early-termination, and infeasibility detection, to be able to efficiently solve LDP relaxations of the form \eqref{eq:ldp-relax}.}
  \label{alg:daqp}
  \begin{algorithmic}[1]
      \Require $M,\underline{d},\overline{d}, \underline{B}, \overline{B}, \lambda_0, \bar{J}$
	\Ensure $u^*, J^*, \mathcal{U}^*, \mathcal{L}^*$
	\algrenewcommand\algorithmicindent{0.825em}%
	\State $\mathcal{U}\leftarrow \underline{B}$, $\mathcal{L} \leftarrow \overline{\mathcal{B}}$, $\lambda \leftarrow \lambda_0$, $\mathcal{E}\leftarrow \underline{\mathcal{B}} \cup \overline{\mathcal{B}}$ 
	\Repeat
	\If{$M_{\mathcal{\mathcal{U}\cup\mathcal{L}}} M_{\mathcal{\mathcal{U}\cup\mathcal{L}}}^T$ is nonsingular}\label{step:nonsingular-start}
	\State $\lambda^* \leftarrow$ solve \eqref{eq:subproblem}
	\If{$\lambda^*_i \geq 0$ $\forall i \in \mathcal{U}$ and $\lambda^*_i \leq 0$ $\forall i \in \mathcal{L}$} \label{step:csp}
	\State $u \leftarrow -M_{\mathcal{U}\cup\mathcal{L}}^T \lambda^*_{\mathcal{U}\cup\mathcal{L}}$
	\If{$\frac{1}{2}\|u\|^2_2 \geq \bar{J}$} infeasible, \Return $\star, \infty, \star, \star$\label{step:early-term}
	\EndIf
	\State $\lambda \leftarrow \lambda^*$
  \State $\mu^+ \leftarrow \overline{d} -M u$,\quad $\mu^-\leftarrow M u - \underline{d}$
	\If{$\mu^+,\mu^- \geq 0$} optimal, \textbf{goto} \ref{step:return}\label{step:primalfeas}
	\Else \quad$j \leftarrow \text{argmin}_{i\notin {\mathcal{U}\cup\mathcal{L}}} \min\{\mu^+_i, \mu^-_i\}$ 
	\If{$\mu^+_j < \mu^-_j$} \:$\mathcal{U} \leftarrow \mathcal{U}\cup\{j\}$
	\EndIf
	\If{$\mu^+_j > \mu^-_j$} \:$\mathcal{L} \leftarrow \mathcal{L}\cup\{j\}$
	\EndIf
	\EndIf
	\Else\: (blocking constraint)
	\State $p \leftarrow \lambda^*-\lambda$\label{step:innerstart} 
    \State$\mathcal{C}\leftarrow\{i\in\mathcal{U}\setminus \mathcal{E}: \lambda^*_i < 0\} \cup \{i\in\mathcal{L}\setminus \mathcal{E}: \lambda^*_i > 0\}$ \label{step:nonsingblocking}
	\State $\lambda, \mathcal{U}, \mathcal{L} \leftarrow$ \textsc{fixComponent}$(\lambda,\mathcal{U},\mathcal{L},\mathcal{C},p)$\label{step:innerstop}
	\EndIf
	\Else $\:(M_{{\mathcal{U}\cup\mathcal{L}}} M_{\mathcal{\mathcal{U}\cup\mathcal{L}}}^T \text{singular})$\label{step:singular-start}
    \State $p \leftarrow$ solve \eqref{eq:subproblem-singular}  \label{step:daqp-sing-dir}
	\State $\mathcal{C} \leftarrow \{i\in \mathcal{U}\setminus \mathcal{E}: p_i < 0\}\cup\{i\in \mathcal{L}\setminus \mathcal{E}: p_i > 0\}$ \label{step:singblocking}
	\If{$\mathcal{C}\setminus \mathcal{E} = \emptyset$} infeasible, \Return $\star, \infty, \star, \star$ 
	\EndIf
	\State $\lambda, \mathcal{U},\mathcal{L} \leftarrow$ \textsc{fixComponent}$(\lambda,\mathcal{U},\mathcal{L},\mathcal{C},p)$\label{step:singular-stop}
	\EndIf
	\Until termination 
	\State \Return $u$, $\frac{1}{2}\|u\|^2_2$, $\mathcal{U}$, $\mathcal{L}$ \label{step:return}
	\hrule\rule{0pt}{1pt}
	\Procedure{fixComponent}{$\lambda, \mathcal{U},\mathcal{L},\mathcal{C},p$}
	\State $j \leftarrow \text{argmin}_{i \in \mathcal{C}} -\lambda_i/p_i $ \label{step:ratio} 
	\State $\lambda \leftarrow \lambda - (\lambda_j/p_j) p$
	\If{$p_j<0$} \:$\mathcal{U} \leftarrow \mathcal{U}\setminus\{j\}$
	\EndIf
	\If{$p_j>0$} \:$\mathcal{L} \leftarrow \mathcal{L}\setminus \{j\}$
	\EndIf
	\State  \Return $\lambda,\mathcal{U},\mathcal{L}$
	\EndProcedure
  \end{algorithmic}
\end{algorithm}
\subsubsection{Double-sided constraints} 
The main focus in \cite{arnstrom2022daqp} was on problems with single-sided constraints of the form ${M u \leq  d}$, while the problem in \eqref{eq:mildp} contains double-sided constraints of the form $\underline{d} \leq M u \leq \bar{d}$. Handling double-sided constraints requires not only keeping track of the components of the dual variable that are free (which are bookkept in a so-called working set $\mathcal{W}$), but also whether the component is allowed to be positive or negative, corresponding to an upper or lower constraints, respectively. Hence, the working set $\mathcal{W}$ is replaced by the working sets $\mathcal{L}$ and $\mathcal{U}$ (where $\mathcal{W}= \mathcal{L}\cup \mathcal{U}$), which contain lower and upper constraints, respectively, that are active (hold with equality). See Section IV.C in \cite{arnstrom2022daqp} for details.

Introducing $\mathcal{L}$ and $\mathcal{U}$ results in the subproblem that is solved in an iteration (Steps \ref{step:csp} and \ref{step:daqp-sing-dir}) being of the form 
\begin{equation}
  \label{eq:subproblem}
  \begin{aligned}
	\lambda^* =  &\:\:\underset{\lambda}{\argmin}&& \frac{1}{2}\lambda^T M M^T \lambda + \overline{d}_{\mathcal{U}}^T \lambda_{\mathcal{U}}+\underline{d}_{\mathcal{L}}^T \lambda_{\mathcal{L}}\\
	 &\text{subject to} &&\lambda_i = 0,\: \forall i\notin \mathcal{L}\cup \mathcal{U},
  \end{aligned}
\end{equation}
when the matrix $M_{\mathcal{U}\cup \mathcal{L}} M_{\mathcal{U}\cup \mathcal{L}}^T$ is non-singular (index sets as subscripts means extracting the corresponding rows from the vector/matrix); and  
\begin{equation}
  \label{eq:subproblem-singular}
  p= 
  \underset{p}{\sol}  
  \left\{
  \begin{aligned}
	M_{\mathcal{U}\cup \mathcal{L}} M_{\mathcal{U}\cup \mathcal{L}}^T \:p_{\mathcal{U}\cup \mathcal{L}}&= 0,\\ 
	\overline{d}_{\mathcal{U}}^T p_{\mathcal{U}} +\underline{d}^T_{\mathcal{L}}p_{\mathcal{L}} < 0,
	\quad p_i=&0\:\: \forall i\notin \mathcal{U}\cup \mathcal{L}
  \end{aligned}
\right\},
\end{equation}
when the matrix $M_{\mathcal{U}\cup \mathcal{L}} M_{\mathcal{U}\cup \mathcal{L}}^T$ is singular. Again, more details about the subproblems are given in \cite{arnstrom2022daqp}.

In DAQP, both \eqref{eq:subproblem} and \eqref{eq:subproblem-singular} are solved efficiently by decomposing $M_{\mathcal{U}\cup \mathcal{L}} M_{\mathcal{U}\cup \mathcal{L}}^T$ with an LDL$^T$ factorization. Moreover, this factorization is recursively updated whenever a constraint is added/removed to/from the working sets $\mathcal{U}$ and $\mathcal{L}$. For details regarding these updates, see Section II.B in \cite{arnstrom2022daqp}.    
\subsubsection{Equality constraints}
To handle equality constraints, which is necessary due to the form of the relaxations in \eqref{eq:ldp-relax}, we make sure that indices corresponding to these, given by the set $\mathcal{E}$, are always contained in the working set (which enforces the constraints to hold with equality.) To this end, we ensure that the set of candidate indices for removal from $\mathcal{U}$ or $\mathcal{L}$, denoted $\mathcal{C}$ (computed at Step \ref{step:nonsingblocking} and \ref{step:singblocking}), never contains any index in $\mathcal{E}$.

\subsubsection{Early-termination and infeasibility detection}
\label{sssec:early-term}
Infeasibility of \eqref{eq:ldp-relax} can be detected whenever $\mathcal{M}_{\mathcal{U}\cup \mathcal{L}} \mathcal{M}_{\mathcal{U}\cup \mathcal{L}}^T$ is singular and $\mathcal{C}= \emptyset$, since then the dual objective function can be made unbounded (see, e.g., Sec. 5.2.2 in \cite{boyd2004convex} for details.)
Moreover, whenever we have dual feasible iterates (i.e., whenever the condition at Step \ref{step:csp} is satisfied) the corresponding primal iterate $u$ will, by duality theory, yield a lower bound $\frac{1}{2} \|u\|^2_2$ of the solution to \eqref{eq:ldp-relax}. As is mentioned in \cite{arnstrom2022daqp} and explored in \cite{fletcher1998numerical}, this lower bound can be used to terminate the solver early if we require the solution to be less than some upper bound $\bar{J}$. If $\frac{1}{2} \|u\|_2^2\geq \bar{J}$, we, hence, consider the problem ``futile'' and early-terminate Algorithm \ref{alg:daqp} at Step \ref{step:early-term}.

\subsubsection{Exact regularization}
\label{sssec:ex-reg}
Nominally, the method in \cite{arnstrom2022daqp} requires the relaxations, and hence the MIQP in \eqref{eq:miqp}, to be strictly convex, i.e., it requires $H \succ 0$. One way of handling singular Hessians is to perform proximal-point iterations \citep{bemporad2018numerically}. This introduces an extra layer of complexity since a sequence of QPs need to be solved for each relaxation. Moreover, the above-mentioned early termination cannot, then, be applied directly since each proximal-point iteration decreases the objective function, resulting in the dual objective function of inner, regularized, QPs not necessarily being a lower bound to the optimal objective function of \eqref{eq:miqp}. 

In hybrid MPC applications, singularity of the Hessian $H$ often originates from binary variables $\delta \in \{0,1\}$ not entering the objective function, called auxiliary variables in MLD systems, since these often encode logical rules that are not directly penalized ($\delta$, hence, only enters the constraints). A naive way of dealing with such singularities is to add a regularizing term $ \epsilon \|\delta\|_2^2$ to the objective. However, this can perturb the true solution, and, more critically in practice, often leads to weakly active constraints in the relaxations, which in turn can lead to numerical instability. We propose instead to add a regularizing term $\epsilon \|\delta-\frac{1}{2}\|_2^2$ to the objective, which does not perturb the solution, and does not encourage weakly active constraints. 
We generalize this by the following proposition.
\begin{prop}[Exact regularization]
    Adding regularizing terms of the form $\|A_i x - \frac{\underline{b}_i + \overline{b}_i}{2}\|_2^2$, $i\in \mathcal{B}$ to the objective function in \eqref{eq:miqp-obj} does not change the solution of \eqref{eq:miqp}.
\end{prop}
\begin{\myproof}
    Denote $q(x)$ the objective function \eqref{eq:miqp-obj} and let 
    \begin{equation}
        \label{eq:reg-problem}
        \begin{aligned}
            \tilde{x}^* =\: &\underset{x}{\argmin} && q(x) + \|A_i x - \tfrac{\underline{b}_i + \overline{b}_i}{2}\|^2_2\\ 
                        &\text{subject to}  && \text{\eqref{eq:miqp-con} and \eqref{eq:miqp-bin}},
        \end{aligned}
    \end{equation}
    
for $i \in \mathcal{B}$. Any $x$ that satisfies the binary constraints \eqref{eq:miqp-bin} also satisfies $\|A_i x- \frac{\underline{b}_i+\overline{b}_i}{2}\|^2_2 = \|\frac{\underline{b}_i - \overline{b}_i}{2}\|_2^2$ for any $i \in \mathcal{B}$.
    Now let $x^*$ be a solution to \eqref{eq:miqp} and assume that $q(x^*) < q(\tilde{x}^*)$. Then, since both $x^*$ and $\tilde{x}^*$ satisfy \eqref{eq:miqp-bin}, we have that 
    \begin{equation*}
        q(x^*) + \|A_i x^* - \tfrac{\underline{b}_i+\overline{b}_i}{2}\|^2_2 <   q(x^*) + \|A_i \tilde{x}^* - \tfrac{\underline{b}_i+\overline{b}_i}{2}\|^2_2,
    \end{equation*}
    which contradicts \eqref{eq:reg-problem}. In conclusion $q(x^*) \geq q(\tilde{x}^*)$, so $\tilde{x}^*$ must be a solution to \eqref{eq:miqp}. 
    \qed
\end{\myproof}
\subsection{Tree exploration}
\label{ssec:branch}
There are three choices that affect the tree exploration in Algorithm \ref{alg:bnb}: 
\textit{node} selection (Step \ref{step:node-selection}), \textit{branch} selection (Step \ref{step:branch-selection}), and \textit{child} selection (Step \ref{step:child-selection}). Specific choices of these, particularly suited for embedded applications, are given below. 

\subsubsection{Node selection} The two most popular search strategies in B\&B are \textit{depth-first} and \textit{best-first}. 
A depth-first search selects the pending node with the highest level $\left|{\overline{\mathcal{B}}}\right|+\left|\underline{\mathcal{B}}\right|$, which often encourages processing of nodes that yield binary feasible solutions. Best-first selects the pending node with the lowest objective function value $\underline{J}$,  which encourages processing of nodes that yield higher quality solutions. While a best-first search often results in fewer nodes being processed in total compared with depth-first, we employ a depth-first because it leads to a reduced memory footprint since the number of pending node is kept low. 
Moreover, since a depth-first search promotes processing children right after their parents have been processed, the inner solver can be \textit{hot-started}, i.e., the working set and the corresponding LDL$^T$ factorization that were formed in the parent can be \textit{directly} reused when solving a child's relaxation, significantly reducing computations. Hot-starting can also be used in best-first, but this requires LDL$^T$ factorizations to be stored for each pending node, leading to a large memory footprint that is, again, not suitable for embedded applications.

Finally, depth-first enables the entire tree to be compactly represented, as is formalized in the following proposition and exemplified below.  
\begin{prop}[Compact tree representation]
    \label{prop:storage}
    If a depth-first search strategy is employed at Step \ref{step:node-selection} in Algorithm \ref{alg:bnb}, a node can be represented with 2 signed integers and in total only $2 n_b + n_b$ signed integers need to be allocated to represent the tree.   
\end{prop}

\begin{pf}
    Let $\tilde{\mathcal{B}}$ be an array containing $n_b$ signed integers. Consider a node at level $\ell$ that was spawned from its parent by fixing $i\in \mathcal{B}$. When this node is processed, let the $\ell$th element of $\tilde{\mathcal{B}}$ be set to index $\pm i$ ($+$ if the upper bound was fixed, $-$ if the lower was fixed). Then when processing any node $(\underline{\mathcal{B}}, \overline{\mathcal{B}})$ at level $\tilde{\ell}$, the first $\tilde{\ell}-1$ elements in $\tilde{\mathcal{B}}$ contains $\underline{\mathcal{B}}$ (the negative elements) and $\overline{\mathcal{B}}$ (the positive elements), since a depth first search ensures that only elements in $\tilde{B}$ at index $\geq \tilde{\ell}$ can have been modified before processing the current node, which means that all singed indices from its parent will be the first $\tilde{\ell}-1$ element of $\tilde{B}$. The only additional information needed to fully retrieve $\underline{\mathcal{B}}$ and $\overline{\mathcal{B}}$, then, is the index that was fixed in the node's parent. 
Hence, the only required storage for a single node is two integers: its level and the (signed) index to add. Moreover, since the number of pending nodes can maximally be $n_b$ when a depth-first search is used, the memory footprint for all pending nodes is maximally $2 n_b$ signed integers. Finally, the buffer $\tilde{\mathcal{B}}$ contains $n_b$ elements, resulting in maximally $2 n_b + n_b$ signed integer being necessary for representing the tree at any iteration.\qed
\end{pf}

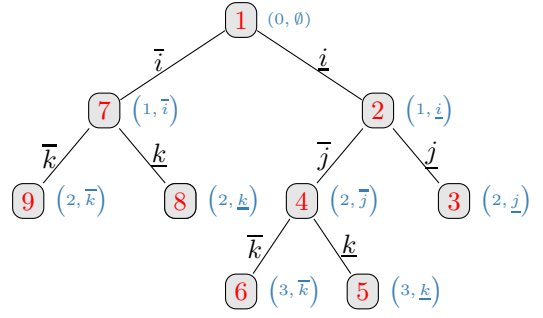
\begin{figure}
  \centering
  \begin{tikzpicture}[level distance=1.5cm,
	level 1/.style={sibling distance=4.5cm},
	level 2/.style={sibling distance=2.5cm},
	level 3/.style={sibling distance=2cm}, 
    every label/.style={set19c2},
    treenode/.style={draw,rounded corners,fill=gray!20,text=red},
    scale=0.8]
    \node[treenode,label=right:{\tiny $\left(0,\emptyset\right)$}] {$1$} 
        child {
            node[treenode,label=right:{\tiny $\left(1,\overline{i}\right)$}] {$7$} 
            child {
                node[treenode,label=right:{\tiny $\left(2,\overline{k}\right)$}] {$9$} 
                edge from parent node [left] {$\overline{k}$}
            }
            child {
                node[treenode] {$8$}
                node[treenode,label=right:{\tiny $\left(2,\underline{k}\right)$}] {$8$} 
                edge from parent node [right] {$\underline{k}$}
            }
            edge from parent node [left,yshift=2pt] {$\overline{i}$}
        }
      child {
          node[treenode,label=right:{\tiny $\left(1,\underline{i}\right)$}] {$2$} 
          child {
              node[treenode,label=right:{\tiny $\left(2,\overline{j}\right)$}] {$4$} 
              child {
                  node[treenode,label=right:{\tiny $\left(3,\overline{k}\right)$}] {$6$} 
                  edge from parent node [left] {$\overline{k}$} 
              }
              child {
                  node[treenode,label=right:{\tiny $\left(3,\underline{k}\right)$}] {$5$} 
                  edge from parent node [right] {$\underline{k}$} 
              }
              edge from parent node [left] {$\overline{j}$} 
          }
          child {
              node[treenode,label=right:{\tiny $\left(2,\underline{j}\right)$}] {$3$} 
              edge from parent node [right] {$\underline{j}$}
          }
          edge from parent node [right,yshift=2pt] {$\underline{i}$}
      };
\end{tikzpicture}
\caption{\small Example of how the branch-and-bound tree can be implicitly represented by two integers per node when a depth-first search is used. The \textcolor{red!70}{number} inside a node corresponds to the order in which it is processed. Indices above edges corresponds to the index that is fixed when moving from the parent to the child. Finally, left of each node is \textcolor{set19c2}{two integers} that are necessary to represent the node, namely, the level and the (signed) index to fix.}
\label{fig:tree}
\end{figure}

To exemplify the result of Proposition \ref{prop:storage}, consider an example with three binary constraints ($n_b =3$) with indices $i$, $j$, and $k$, which requires $\tilde{\mathcal{B}}$ to contain at the maximum three elements, i.e.,  
$\tilde{\mathcal{B}} = [\star, \star, \star]$. A possible B\&B-tree for this scenario is shown in Figure \ref{fig:tree}, and how this tree can be compactly represented will now be motivated.

\newcommand{\RNum}[1]{\uppercase\expandafter{\romannumeral #1\relax}}

Consider the case when node $5$ should be processed. Then $\tilde{\mathcal{B}}  = \{\underline{i}, \overline{j}, \star\}$ since node $4$ was processed right before. The only information required at node $5$ is its placement in the tree (level 3) and which index should be fixed (index $\underline{k}$). All other information is implicitly stored in $\tilde{\mathcal{B}}$ since we know that all elements up to level 3 ($\underline{i}$ and $\overline{j}$) will be correct (only nodes at the same level or below the current node could have been processed in between the current node's parent and itself since we use a depth-first search).

\subsubsection{Branch selection}
We use a lexicographic selection rule that picks the smallest index in $\mathcal{B}$ that corresponds to a constraint that is satisfied; that is, the branching index $i$ is selected as $i =  \min_j \left\{j\in \mathcal{B}: M_j u \notin \{\underline{d}_j, \bar{d}_j\} \right\}$. This allows the user to order the constraints according to their branching priority.
There are several branching rules that are more advanced, for example \textit{strong branching}, \textit{reliability branching}, or \textit{hybrid branching} \citep{achterberg2005branching}. While the simpler lexicographical rule might lead to more processed nodes compared with the more complex rules mentioned above, we use it for its simplicity and, hence, its small overhead, particularly suitable for embedded applications.
Also, in the context of MPC of hybrid systems, a simple lexicographic selection rule with an ordering according to a time index can be very effective \citep{bemporad1999efficient}.
\subsubsection{Child selection} We explore the child that corresponds to the bound that is closest to being satisfied in the parent. That is, if $\tilde{u}$ is the solution in the parent $(\underline{\mathcal{B}}, \overline{\mathcal{B}})$ and constraint $i$ is selected to be branched upon, the node $(\underline{\mathcal{B}}\cup\{i\}, \overline{\mathcal{B}})$ is processed first if $M_i \tilde{u} \leq \frac{\underline{d}_i+\bar{d}_i}{2}$; otherwise, the node $(\underline{\mathcal{B}}, \overline{\mathcal{B}}\cup\{i\})$ is processed first. 
\begin{\myremark}[Most fractional branching rule]
    A similar concept to the one used for child selection can be used in branch selection (known as the ``most fractional rule'' \cite{achterberg2005branching}). However, this requires the product $M_i \tilde{u}$ to be computed \textit{for all} branching candidates (while only a single product needs to be computed when it is just used for child selection). Moreover, as is shown in \cite{achterberg2005branching}, the most fractional rule for node selection seldom improves the number of processed nodes compared with selecting the branching index randomly. 
\end{\myremark}

\subsection{Complete algorithm}
\label{ssec:complete-method}
We are now ready to present the main contribution of this paper.
Algorithm \ref{alg:bnb-full} concretizes Algorithm \ref{alg:bnb} by adding the details described in the preceding subsections. Namely, by solving the relaxations using Algorithm \ref{alg:daqp} and by using the search strategy outlined in Section \ref{ssec:branch}.
In addition, the algorithm contains a step for transforming an MIQP of the form \eqref{eq:miqp} into an MILDP of the form \eqref{eq:mildp} at Step \ref{step:transform}, and a step for retrieving the solution to \eqref{eq:miqp} from a solution of \eqref{eq:mildp} at Step \ref{step:retreive}. 
\begin{algorithm}[H]
  \algrenewcommand\algorithmicindent{1.2em}
  \caption{B\&B method for solving the MIQP in \eqref{eq:miqp}}
  \label{alg:bnb-full}
  \begin{algorithmic}[1]
	\Require $R,f,A,\underline{b},\bar{b}, \mathcal{B}$
    \Ensure $x^*$
    \State $M \leftarrow A R^{-1}$, $v\leftarrow R^{-T}f$,  $\underline{d} \leftarrow \underline{b}+Mv$, $\overline{d} \leftarrow \overline{b} + M v$ \label{step:transform} 
	\State $\bar{u}\leftarrow 0$, $\bar{J} \leftarrow \infty$, $\mathcal{T} \leftarrow \{(\emptyset, \emptyset)$\}
	\While{$\mathcal{T}\neq \emptyset$}
	\State $(\underline{\mathcal{B}},\overline{\mathcal{B}}) \leftarrow$ pop from $\mathcal{T}$
    \State $u, J, \mathcal{U},\mathcal{L} \leftarrow$ solve \eqref{eq:ldp-relax} using Algorithm \ref{alg:daqp}
	\If{$J \geq  \bar{J}$} \textbf{continue}
	\EndIf
	\If{$\mathcal{B} \subseteq (\mathcal{U}\cup \mathcal{L})$} 
	\State $\bar{u} \leftarrow u, \bar{J} \leftarrow J$
	\Else \Comment branching
  \State $i\leftarrow \min \{ i\in \mathcal{B}: i\notin \mathcal{U}\cup \mathcal{L}\}$  
	\If{$M_i u \leq \frac{\underline{d}_i+\bar{d}_i}{2}$}
	\State push $(\underline{\mathcal{B}},\overline{\mathcal{B}}\cup\{i\})$ to $\mathcal{T}$; push $(\underline{\mathcal{B}}\cup\{i\},\overline{\mathcal{B}})$ to $\mathcal{T}$
	\Else
	\State push $(\underline{\mathcal{B}}\cup\{i\},\overline{\mathcal{B}})$ to $\mathcal{T}$; push $(\underline{\mathcal{B}},\overline{\mathcal{B}}\cup\{i\})$ to $\mathcal{T}$
	\EndIf
	\EndIf
	\EndWhile
    \State \Return $x^* \leftarrow - R^{-1}(\bar{u}-v)$\label{step:retreive}
  \end{algorithmic}
\end{algorithm}

\subsection{Complexity certification}
When considering MIQPs originating from hybrid \textit{linear} MPC, $f,\overline{b},$ and $\underline{b}$ are affine functions of a parameter $\theta \in \mathbb{R}^p$. This structure can be exploited by complexity-certification methods to determine \emph{tight} worst-case guarantees on the number of required computations Algorithm \ref{alg:bnb-full} requires for solving \emph{any} MIQP generated by $\theta$ \citep{shoja2022overall}.
The complexity-certification in \cite{shoja2022overall} can be directly applied to Algorithm \ref{alg:bnb-full} since it fulfills the requirements specified therein. One such requirement is fulfilled by the inner solver DAQP being certifiable by the framework in \cite{arnstrom2022unifying}.
Detailed parametric complexity-certification results provided by the method in \cite{shoja2022overall} for the proposed solver is out of scope here and will be considered in future work.

\section{Numerical Experiments}
In this section we report results from numerical experiments for a C implementation\footnote{available at \url{https://github.com/darnstrom/daqp}}
of Algorithm \ref{alg:bnb-full} (BnB-DAQP). In Section \ref{ssec:result-rand} we compare BnB-DAQP with the state-of-the-art solver Gurobi on a set of randomly generated MIQPs. Then, in Section \ref{ssec:result-embed}, we highlight the embeddability of BnB-DAQP by implementing it on an MCU with limited memory and computational resources, which is then used for hybrid MPC. 
Code for all experiments is available online at \url{https://github.com/darnstrom/ifac2023-bnbdaqp}.
\begin{rem}[Comparing against additional MIQP solvers]
    We only compare against Gurobi because most MIQP solvers in the MPC literature \citep[e.g.,][]{bemporad2015solving,bemporad2018numerically,hespanhol2019structure,liang2020early,stellato2018embedded} are not publicly accessible or is not readily embeddable. Hence, the open-source C implementation\footnotemark[1] of our solver is a major contribution of this work. Nevertheless, most of the above-mentioned solvers are compared with Gurobi; so a rough comparison could be made implicitly by comparing their relative performance with the results in Section \ref{ssec:result-rand}.
\end{rem}
\subsection{Random mixed-integer QPs}
\label{ssec:result-rand}
First, we compare BnB-DAQP 
with the state-of-the-art commercial solver Gurobi (version 9.5.2) on a set of randomly generated MIQPs. The problems are of the form
\begin{subequations}
  \label{eq:miqp-random}
  \begin{align}
   &\underset{x}{\text{minimize}}&&\frac{1}{2} x^T H x + f^T x \label{eq:miqprandom-obj}\\
   &\text{subject to} &&\underline{b} \leq A x \leq \bar{b} \label{eq:miqpradnom-con}\\
   & && x_i \in \left\{0, 1  \right\},\quad \forall i \in \{1,\dots, n_b\},
  \end{align}
\end{subequations}
where elements of $A,\bar{b}$, and $\underline{b}$ are generated as  
$A\sim\mathcal{N}(0,1)$, $\bar{b}\sim \mathcal{U}(0,20)$, and $\underline{b}\sim \mathcal{U}(-20,0)$, with $\mathcal{N}(\mu,\sigma)$ denoting a normal distribution with mean $\mu$ and standard deviation $\sigma$, and $\mathcal{U}(l,u)$ denotes the uniform distribution over the interval $[l,u]$. Moreover, the Hessian $H$ was generated with the MATLAB function \texttt{sprandsym} with density $1$ and condition number $10^{-4}$, and the linear term $f$ was partitioned as $f = \left(\begin{smallmatrix}
 f_b \\ f_c 
\end{smallmatrix}\right)$, with  $f_b \sim -|\mathcal{N}(0,10^2)|$ and $f_c \sim \mathcal{N}(0,10^2)$; the negativity of $f_b$ was enforced to counteract a bias of $x^*_i=0$,  $i\in \{1,\dots, n_b\}$, since $H \succ 0$. 

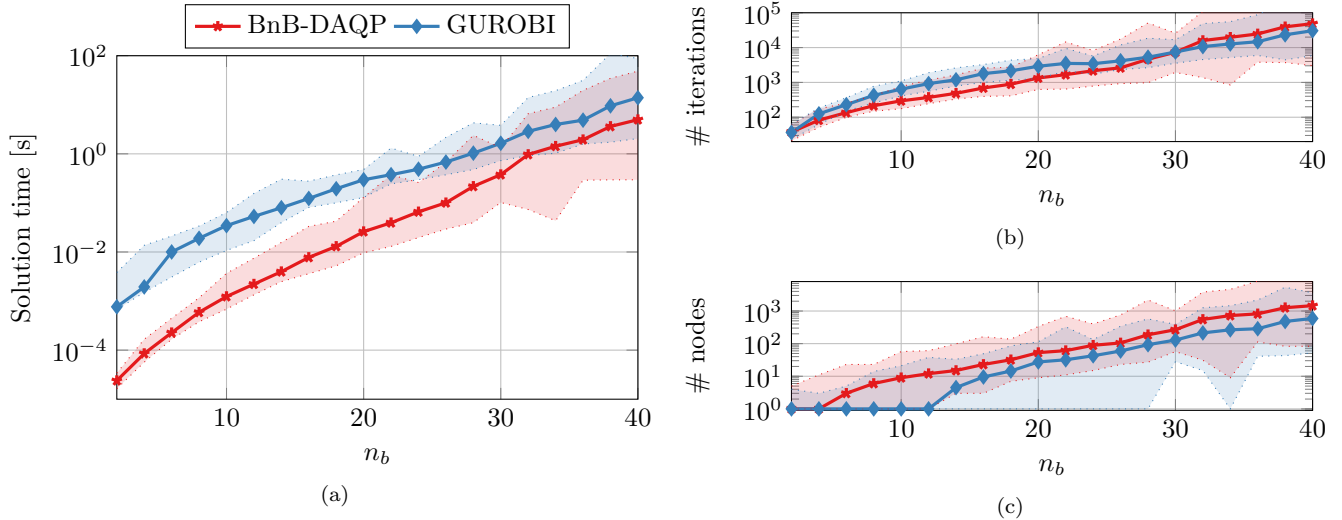
\begin{figure*}[htpb] 
  \centering
  \begin{minipage}{.45\linewidth}
  \subfloat[\label{subfig:random-time}]{%
	\begin{tikzpicture}[scale=1]
	  \begin{axis}[xmin=2,xmax=40,
		ymode=log,
		ymin=10^-5,ymax=10^2,
		xlabel={$n_b$},
		ylabel={Solution time [s]},
		legend style={at ={(0.5,1.15)},anchor=north}, ymajorgrids,xmajorgrids,
		x post scale=1,
		y post scale=0.8,
		legend cell align={left},legend columns=2,
		]
	 	\addplot [set19c1,very thick, mark=star] table [x={nb}, y={daqpmed}] {\randomTimes}; 
		\addplot [set19c2,very thick, mark=diamond*] table [x={nb}, y={grbmed}] {\randomTimes};
	 	\addplot [set19c1,dotted, name path=updaqp] table [x={nb}, y={daqpwc}] {\randomTimes}; 
	 	\addplot [set19c1,dotted, name path=lowdaqp] table [x={nb}, y={daqpbc}] {\randomTimes}; 
	 	\addplot [set19c2,dotted, name path=upgrb] table [x={nb}, y={grbwc}] {\randomTimes}; 
	 	\addplot [set19c2,dotted, name path=lowgrb] table [x={nb}, y={grbbc}] {\randomTimes}; 
		\addplot [set19c1,fill opacity=0.15] fill between [of=updaqp and lowdaqp];
		\addplot [set19c2,fill opacity=0.15] fill between [of=upgrb and lowgrb];
		\legend{BnB-DAQP, GUROBI}
	  \end{axis}
	\end{tikzpicture}
  }
\end{minipage}\qquad
  \begin{minipage}{.45\linewidth}
  \subfloat[\label{subfig:random-iter}]{%
	\begin{tikzpicture}[scale=1]
	  \begin{axis}[xmin=2,xmax=40,
		scale=0.7,
		ymode=log,
		ymin=2*10^1,ymax=10^5,
		xlabel={$n_b$},
		ylabel={\# iterations},
		legend style={at ={(1,0)},anchor=south east}, ymajorgrids,xmajorgrids,
		x post scale=1,
		y post scale=0.3,
		legend cell align={left},legend columns=3,
		legend style={nodes={scale=0.65, transform shape}},
		]
	 	\addplot [set19c1,very thick, mark=star] table [x={nb}, y={daqpmed}] {\randomIters}; 
		\addplot [set19c2,very thick, mark=diamond*] table [x={nb}, y={grbmed}] {\randomIters};
	 	\addplot [set19c1,dotted, name path=updaqp] table [x={nb}, y={daqpwc}] {\randomIters}; 
	 	\addplot [set19c1,dotted, name path=lowdaqp] table [x={nb}, y={daqpbc}] {\randomIters}; 
	 	\addplot [set19c2,dotted, name path=upgrb] table [x={nb}, y={grbwc}] {\randomIters}; 
	 	\addplot [set19c2,dotted, name path=lowgrb] table [x={nb}, y={grbbc}] {\randomIters}; 
		\addplot [set19c1,fill opacity=0.15] fill between [of=updaqp and lowdaqp];
		\addplot [set19c2,fill opacity=0.15] fill between [of=upgrb and lowgrb];
	  \end{axis}
	\end{tikzpicture}
} \\
\subfloat[\label{subfig:random-nodes}]{%
	\begin{tikzpicture}[scale=1]
	  \begin{axis}[xmin=2,xmax=40,
		scale=0.7,
		ymode=log,
		ymin=0.9*10^0,ymax=0.8*10^4,
		xlabel={$n_b$},
		ylabel={\# nodes},
		legend style={at ={(1,0)},anchor=south east}, ymajorgrids,xmajorgrids,
		x post scale=1,
		y post scale=0.3,
		legend cell align={left},legend columns=3,
		legend style={nodes={scale=0.65, transform shape}},
		]
	 	\addplot [set19c1,very thick, mark=star] table [x={nb}, y={daqpmed}] {\randomNodes}; 
		\addplot [set19c2,very thick, mark=diamond*] table [x={nb}, y={grbmed}] {\randomNodes};
	 	\addplot [set19c1,dotted, name path=updaqp] table [x={nb}, y={daqpwc}] {\randomNodes}; 
	 	\addplot [set19c1,dotted, name path=lowdaqp] table [x={nb}, y={daqpbc}] {\randomNodes}; 
	 	\addplot [set19c2,dotted, name path=upgrb] table [x={nb}, y={grbwc}] {\randomNodes}; 
	 	\addplot [set19c2,dotted, name path=lowgrb] table [x={nb}, y={grbbc}] {\randomNodes}; 
		\addplot [set19c1,fill opacity=0.15] fill between [of=updaqp and lowdaqp];
		\addplot [set19c2,fill opacity=0.15] fill between [of=upgrb and lowgrb];
	  \end{axis}
	\end{tikzpicture}
}
\end{minipage}
\caption{Solution time, number of solved relaxtions, and total number of inner iterations when solving randomly generated MIQPs of the form \eqref{eq:miqp-random}. The number of decision variables is $n=5 n_b$ and the number of (double-sided) constraints is $m= 10 n_b$ (for example, the largest problem has dimension $n_b=40$, $n=200$, $m=400$). For each value of $n_b$, 50 MIQP instances were generated and solved. Solid lines mark the median and dotted lines mark the best- and worst-case results for these 50 solves. Both solvers were exectued on an Intel 2.7 GHz i7-7500U CPU.}
\label{fig:random-fig} 
\end{figure*}

The result for solving MIQPs of the form \eqref{eq:miqp-random} for varying problem dimensions is shown in Figure \ref{fig:random-fig}.
Figure \ref{subfig:random-time} shows that BnB-DAQP outperforms Gurobi on the considered problems/dimensions (representable of problems commonly encountered in hybrid MPC applications), both concerning average solution time and worst-case solution time. From Figure \ref{subfig:random-nodes} we see that the tree exploration in Gurobi is more effective, as is to be expected since the tree exploration in BnB-DAQP favours simplicity. Despite this, however, the number of iterations (the number of solved linear equation systems), shown in Figure \ref{subfig:random-iter}, is initially lower for BnB-DAQP.

\subsection{Embedded MPC}
\label{ssec:result-embed}
To illustrate the proposed solver's embeddability, we apply it to a hybrid MPC example, with the controller running on an MCU with limited computing power and memory. Specifically, we use an STM32F411 MCU, running at 84 MHz and with 512 kB Flash memory and 128 kB of SRAM. The MCU has no data caches and an FPU that only supports floating point operations in \textit{single} precision.

\begin{figure}[H]
    \centering
    \begin{tikzpicture}[scale=0.75, transform shape]
        \tikzstyle{damper}=[decoration={markings,
  mark connection node=dmp,
  mark=at position 0.5 with
  {
    \node (dmp) [inner sep=0pt,transform shape,rotate=-90,minimum width=6pt,minimum height=1pt,draw=none] {};
    \draw [] ($(dmp.north east)+(2pt,0)$) -- (dmp.south east) -- (dmp.south west) -- ($(dmp.north west)+(2pt,0)$);
    \draw [] ($(dmp.north)+(0,-2pt)$) -- ($(dmp.north)+(0,2pt)$);
  }
}, decorate]
        \coordinate (hitch) at (0.81,1);
        \coordinate (anglehelp) at (0.81,3.25);
        \coordinate (pendulumtip) at (-0.5,3);
        \coordinate (wallup) at (-2,3.5);
        \coordinate (wallmiddle) at (-2,2.75);
        \coordinate (walldown) at (-2,2);
        \draw[fill=black] (0.81,1) circle (3.5pt);
        \draw[line width=2pt] (hitch) -- (pendulumtip);
        \draw[dashed] (hitch) -- (anglehelp);
        \draw[thick] (wallup) -- (walldown);
        \draw[red,-latex,line width=1pt] (wallmiddle) -- node[above]{$F_\text{con}$}++(0.75,0);
        \draw[<->] (-.55,0.61) -- node[above]{$d$}(-2,0.61) ;
        \draw[fill=gray!10] (0,0.25) -- (1.62,0.25)--(1.62,1)--(0,1)--(0,0.25); 
        \draw[fill=black] (pendulumtip) circle (3.5pt);
        \draw[red,-latex,line width=1pt] (1.5,0.61) -- node[above]{$F$}(2.5,0.61);
        \draw[|-latex,set19c2] (-.5,.61) -- node[above]{$x$}(0.81,.61);
        \pic [set19c2,draw,-latex, angle radius=11mm, angle eccentricity=1.2] {angle = anglehelp--hitch--pendulumtip};
    \draw
        [
        decoration={
            coil,
            segment length = 0.5mm,
            amplitude = 1mm,
            aspect = 0.25,
            post length = 1mm,
        pre length = 1mm},
        decorate] (-2,3) -- ++(-.75,0);
    \draw [damper](-2.75,2.5) -- (-2,2.5); 
    \node[set19c2] (test) at (0.5,2.3) {$\phi$};
    \draw[fill=white,thick] (0.35,0.2) circle (4.5pt);
    \draw[fill=white,thick] (1.25,0.2) circle (4.5pt);
    \fill [pattern = north west lines] (-2.5,-.2) rectangle ++(5,.2);
    \draw[thick] (-2.5,0) -- ++(5,0);
    \fill [pattern = north west lines] (-2.75,3.5) rectangle ++(-0.2,-2.5);
    \draw[thick] (-2.75,3.5) -- ++(0,-2.5);
\end{tikzpicture}
\caption{Inverted pendulum on a cart with contact forces.}
\label{fig:invpend-tikz}
\end{figure}
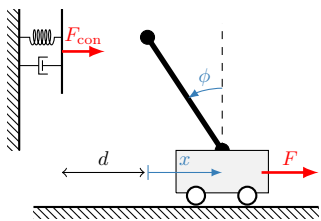

We consider the hybrid MPC of a linearized inverted pendulum on a cart, surrounded by a wall (giving rise to contact forces), considered in, for example, \cite{marcucci2020warm}, and visualized in Figure \ref{fig:invpend-tikz}. The control goal is to stabilize the pendulum in the upright position ($\phi=0$) at the origin ($x=0$). The control is a force $F$ directly applied to the cart. There is also a contact force $F_{\text{con}}$ present when the tip of the pendulum touches the wall, necessitating binary variables in the model.

The specific control task considered in the experiment is a recovery task, where the pendulum is initialized on a collision course ($x = 0, \dot{x} = -1$) with the wall. Moreover, after two seconds another impulse to the velocity of the cart ($\dot{x} = -0.7$) is applied, setting it on yet another collision course with the wall. A similar scenario was considered in \cite{marcucci2020warm}. 

An MPC with a prediction/control horizon 6, where each time step is $0.1$ seconds, was used to control the platform. We impose the control constraint $|F| \leq 1$ and the state constraints $|x| \leq d$ , $|\dot{x}|\leq 1$, $|\phi| \leq \frac{\pi}{10}$, $|\dot{\phi}| \leq 1$ at each time step. Additionally, there are several constraints that contain binary variables to model the contact force. The dimensions of the resulting MIQP problems are $n= 24$, $m=96$, and $n_b = 12$. 
The same weights and model parameters as in \cite{marcucci2020warm} were used. More details about the MPC problem, for example how the constraints from the contact force are formulated in the MIQP, are also given in \cite{marcucci2020warm}.

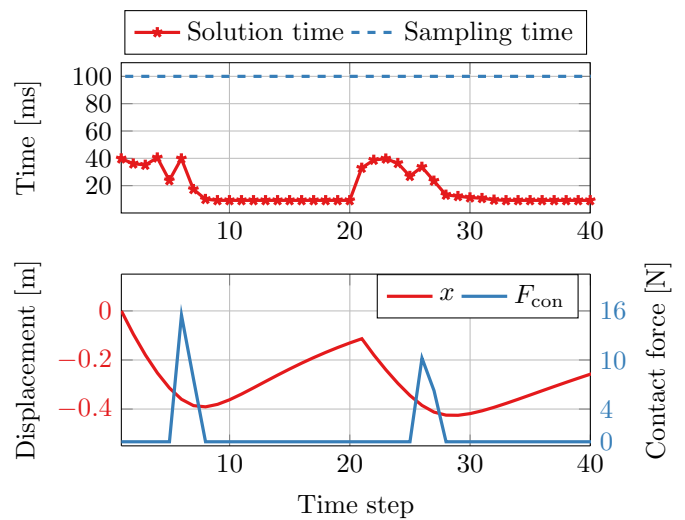
\begin{figure}
  \centering
	\begin{tikzpicture}[scale=1]
	  \begin{axis}[xmin=1,xmax=40,
		ymin=0,ymax=110,
        ytick={20,40,60,80,100},
        xlabel={},
		ylabel={Time [ms]},
		legend style={at ={(0.5,1.35)},anchor=north}, ymajorgrids,yminorgrids,xmajorgrids,
		x post scale=0.9,
		y post scale=0.35,
		legend cell align={left},legend columns=2,
		]
	 	\addplot [set19c1,very thick, mark=star] table [x={step}, y={time}] {\invpend}; 
        \addplot[mark=none,very thick, dashed, set19c2,domain=0:100] {100};
		\legend{Solution time, Sampling time}
	  \end{axis}
      \node (pad) at (10,0) {};
	\end{tikzpicture}
	\begin{tikzpicture}[scale=1]
        \pgfplotsset{set layers}
	  \begin{axis}[xmin=1,xmax=40,
		ymin=-0.55,ymax=0.15,
		xlabel={Time step},
		legend style={at ={(0.975,1)},anchor=north east},ymajorgrids ,xmajorgrids,
		x post scale=1,
        axis y line*=left,
        ylabel={Displacement [m]},
		x post scale=0.9,
        ytick={-0.4,-0.2,0},
		y post scale=0.4,
        y tick label style={color=set19c1},
		legend cell align={left},legend columns=2,
		]
	 	\addplot [set19c1,very thick] table [x={step}, y={x1}] {\invpend}; 
        \addplot [set19c2,very thick] {-10}; 
        \legend{$x$,$F_{\text{con}}$}
	  \end{axis}
	  \begin{axis}[xmin=1,xmax=40,
		ymin=-0.5,ymax=20.5,
        axis y line*=right,
        axis x line=none,
		xlabel={Time step},
        ytick={0,4,10,16},
        ylabel={Contact force [N]},
        ylabel near ticks,
		legend style={at ={(0.5,1.15)},anchor=north}, xmajorgrids,
		x post scale=0.9,
		y post scale=0.4,
        y tick label style={color=set19c2},
		legend cell align={left},legend columns=2,
		]
	 	\addplot [set19c2,very thick] table [x={step}, y={u2}] {\invpend}; 
	  \end{axis}
	\end{tikzpicture}
    \caption{Recovery task of a inverted pendulum (starting with $x=0$ and $\dot{x} = -1$) and an additional impulse in the velocity at time step $20$ with $\dot{x} = -0.7$. The dimensions of the MIQPs are $n=24$, $m=96$ and $n_b =12$, and the proposed MIQP solver was executed on an STM32F411 MCU.}
\label{fig:invpend-res} 
\end{figure}

The solution times for computing a control in the simulated scenario is shown in Figure \ref{fig:invpend-res}, where it can be seen that the controller is able to operate well below the sampling time of 0.1 seconds (i.e., a control frequency of 10 Hz). Note that the solution time also includes forming the problem from a given state estimate. Figure \ref{fig:random-fig} also shows that there are a non-zero contact force at some time steps, meaning that the MIQPs are non-trivial. The total memory footprint (including problem data and code for logging) used on the MCU was 25.4 kB.

\section{Conclusion}
We have proposed a mixed-integer QP solver that is suitable for use in embedded applications such as hybrid model predictive control (MPC). The solver is based on branch-and-bound, and the recently proposed dual active-set solver DAQP \citep{arnstrom2022unifying} is used for solving QP relaxations. We showed that the proposed solver outperforms Gurobi on small to medium-sized MIQPs, commonly encountered in embedded hybrid MPC applications. Finally, the embeddability of the solver was shown by successfully using it on an MCU with limited memory and computing power for the MPC of an inverted pendulum on a cart with contact forces. 
A {C implementation} of the proposed solver is available for download at {\url{https://github.com/darnstrom/daqp}}.
\linespread{1.0}\selectfont
\bibliography{lib}

\begin{thebibliography}{25}
\providecommand{\natexlab}[1]{#1}
\providecommand{\url}[1]{\texttt{#1}}
\providecommand{\urlprefix}{URL }
\expandafter\ifx\csname urlstyle\endcsname\relax
  \providecommand{\doi}[1]{doi:\discretionary{}{}{}#1}\else
  \providecommand{\doi}{doi:\discretionary{}{}{}\begingroup
  \urlstyle{rm}\Url}\fi

\bibitem[{Achterberg et~al.(2005)Achterberg, Koch, and
  Martin}]{achterberg2005branching}
Achterberg, T., Koch, T., and Martin, A. (2005).
\newblock Branching rules revisited.
\newblock \emph{Operations Research Letters}, 33(1), 42--54.

\bibitem[{Arnström and Axehill(2022)}]{arnstrom2022unifying}
Arnström, D. and Axehill, D. (2022).
\newblock A unifying complexity certification framework for active-set methods
  for convex quadratic programming.
\newblock \emph{IEEE Transactions on Automatic Control}, 67(6), 2758--2770.
\newblock \doi{10.1109/TAC.2021.3090749}.

\bibitem[{Arnström et~al.(2022)Arnström, Bemporad, and
  Axehill}]{arnstrom2022daqp}
Arnström, D., Bemporad, A., and Axehill, D. (2022).
\newblock A dual active-set solver for embedded quadratic programming using
  recursive {LDL}$^{T}$ updates.
\newblock \emph{IEEE Transactions on Automatic Control}, 67(8), 4362--4369.
\newblock \doi{10.1109/TAC.2022.3176430}.

\bibitem[{Axehill and Hansson(2006)}]{axehill2006mixed}
Axehill, D. and Hansson, A. (2006).
\newblock A mixed integer dual quadratic programming algorithm tailored for
  {MPC}.
\newblock In \emph{Proceedings of the 45th IEEE Conference on Decision and
  Control}, 5693--5698. IEEE.

\bibitem[{Axehill and Hansson(2008)}]{axehill2008dual}
Axehill, D. and Hansson, A. (2008).
\newblock A dual gradient projection quadratic programming algorithm tailored
  for model predictive control.
\newblock In \emph{2008 47th IEEE Conference on Decision and Control},
  3057--3064. IEEE.

\bibitem[{Axehill and Morari(2010)}]{axehill2010improved}
Axehill, D. and Morari, M. (2010).
\newblock Improved complexity analysis of branch and bound for hybrid {MPC}.
\newblock In \emph{49th IEEE Conference on Decision and Control (CDC)},
  4216--4222. IEEE.

\bibitem[{Bemporad(2015)}]{bemporad2015solving}
Bemporad, A. (2015).
\newblock Solving mixed-integer quadratic programs via nonnegative least
  squares.
\newblock \emph{IFAC-PapersOnLine}, 48(23), 73--79.

\bibitem[{Bemporad et~al.(1999)Bemporad, Mignone, and
  Morari}]{bemporad1999efficient}
Bemporad, A., Mignone, D., and Morari, M. (1999).
\newblock An efficient branch and bound algorithm for state estimation and
  control of hybrid systems.
\newblock In \emph{1999 European Control Conference (ECC)}, 557--562. IEEE.

\bibitem[{Bemporad and Morari(1999)}]{bemporad1999control}
Bemporad, A. and Morari, M. (1999).
\newblock Control of systems integrating logic, dynamics, and constraints.
\newblock \emph{Automatica}, 35(3), 407--427.

\bibitem[{Bemporad and Naik(2018)}]{bemporad2018numerically}
Bemporad, A. and Naik, V.V. (2018).
\newblock A numerically robust mixed-integer quadratic programming solver for
  embedded hybrid model predictive control.
\newblock \emph{IFAC-PapersOnLine}, 51(20), 412--417.

\bibitem[{Bertsimas and Stellato(2022)}]{bertsimas2022online}
Bertsimas, D. and Stellato, B. (2022).
\newblock Online mixed-integer optimization in milliseconds.
\newblock \emph{INFORMS Journal on Computing}.

\bibitem[{Boyd and Vandenberghe(2004)}]{boyd2004convex}
Boyd, S. and Vandenberghe, L. (2004).
\newblock \emph{Convex optimization}.
\newblock Cambridge university press.

\bibitem[{Ferreau et~al.(2014)Ferreau, Kirches, Potschka, Bock, and
  Diehl}]{ferreau2014qpoases}
Ferreau, H.J., Kirches, C., Potschka, A., Bock, H.G., and Diehl, M. (2014).
\newblock {qpOASES}: A parametric active-set algorithm for quadratic
  programming.
\newblock \emph{Mathematical Programming Computation}, 6(4), 327--363.

\bibitem[{Fletcher and Leyffer(1998)}]{fletcher1998numerical}
Fletcher, R. and Leyffer, S. (1998).
\newblock Numerical experience with lower bounds for {MIQP} branch-and-bound.
\newblock \emph{SIAM Journal on Optimization}, 8(2), 604--616.

\bibitem[{Frick et~al.(2015)Frick, Domahidi, and Morari}]{frick2015embedded}
Frick, D., Domahidi, A., and Morari, M. (2015).
\newblock Embedded optimization for mixed logical dynamical systems.
\newblock \emph{Computers \& Chemical Engineering}, 72, 21--33.

\bibitem[{Frison and Diehl(2020)}]{frison2020hpipm}
Frison, G. and Diehl, M. (2020).
\newblock {HPIPM}: a high-performance quadratic programming framework for model
  predictive control.
\newblock \emph{IFAC-PapersOnLine}, 53(2), 6563--6569.

\bibitem[{Hespanhol et~al.(2019)Hespanhol, Quirynen, and
  Di~Cairano}]{hespanhol2019structure}
Hespanhol, P., Quirynen, R., and Di~Cairano, S. (2019).
\newblock A structure exploiting branch-and-bound algorithm for mixed-integer
  model predictive control.
\newblock In \emph{2019 18th European Control Conference (ECC)}, 2763--2768.
  IEEE.

\bibitem[{Land and Doig(1960)}]{land1960automatic}
Land, A. and Doig, A. (1960).
\newblock An automatic method of solving discrete programming problems.
\newblock \emph{Econometrica: Journal of the Econometric Society}, 497--520.

\bibitem[{Liang et~al.(2020)Liang, Di~Cairano, and Quirynen}]{liang2020early}
Liang, J., Di~Cairano, S., and Quirynen, R. (2020).
\newblock Early termination of convex {QP} solvers in mixed-integer programming
  for real-time decision making.
\newblock \emph{IEEE Control Systems Letters}, 5(4), 1417--1422.

\bibitem[{Marcucci and Tedrake(2020)}]{marcucci2020warm}
Marcucci, T. and Tedrake, R. (2020).
\newblock Warm start of mixed-integer programs for model predictive control of
  hybrid systems.
\newblock \emph{IEEE Transactions on Automatic Control}, 66(6), 2433--2448.

\bibitem[{Naik and Bemporad(2017)}]{naik2017embedded}
Naik, V.V. and Bemporad, A. (2017).
\newblock Embedded mixed-integer quadratic optimization using accelerated dual
  gradient projection.
\newblock \emph{IFAC-PapersOnLine}, 50(1), 10723--10728.

\bibitem[{Patrinos and Bemporad(2013)}]{patrinos2013accelerated}
Patrinos, P. and Bemporad, A. (2013).
\newblock An accelerated dual gradient-projection algorithm for embedded linear
  model predictive control.
\newblock \emph{IEEE Transactions on Automatic Control}, 59(1), 18--33.

\bibitem[{Shoja et~al.(2022)Shoja, Arnström, and Axehill}]{shoja2022overall}
Shoja, S., Arnström, D., and Axehill, D. (2022).
\newblock Overall complexity certification of a standard branch and bound
  method for mixed-integer quadratic programming.
\newblock In \emph{2022 American Control Conference (ACC)}, 4957--4964.
\newblock \doi{10.23919/ACC53348.2022.9867176}.

\bibitem[{Stellato et~al.(2018)Stellato, Naik, Bemporad, Goulart, and
  Boyd}]{stellato2018embedded}
Stellato, B., Naik, V.V., Bemporad, A., Goulart, P., and Boyd, S. (2018).
\newblock Embedded mixed-integer quadratic optimization using the {OSQP}
  solver.
\newblock In \emph{2018 European Control Conference (ECC)}, 1536--1541. IEEE.

\bibitem[{Takapoui et~al.(2020)Takapoui, Moehle, Boyd, and
  Bemporad}]{takapoui2020simple}
Takapoui, R., Moehle, N., Boyd, S., and Bemporad, A. (2020).
\newblock A simple effective heuristic for embedded mixed-integer quadratic
  programming.
\newblock \emph{International journal of control}, 93(1), 2--12.

\end{thebibliography}

\end{document}